\documentclass[reqno,12pt]{amsart}
\usepackage{latexsym,amsmath,amsthm,amssymb,amscd}
\usepackage{hyperref}

\newcommand{\leftexp}[2]{{\vphantom{#2}}^{#1}{#2}}	
\renewcommand{\to}{\rightarrow}
\newcommand{\F}{\ensuremath{\mathsf{F}}} 
\newcommand{\M}{\ensuremath{\mathsf{M}}} 
\renewcommand{\L}{\ensuremath{\mathsf{L}}} 
\newcommand{\R}{\ensuremath{\mathsf{R}}} 
\newcommand{\W}{\ensuremath{\mathsf{W}}} 
\newcommand{\C}{\ensuremath{\mathcal{C}}} 
\newcommand{\sfC}{\ensuremath{\mathsf{C}}} 
\newcommand{\T}{\ensuremath{\mathbb{T}}} 
 
\usepackage[all]{xy}
\xyoption{2cell}
\xyoption{curve}
\UseTwocells
\input{diagxy}

\newcommand{\A}{A_\bullet}

\newcommand{\G}{\mathbb{G}}
\newcommand{\N}{\mathbb{N}}

\renewcommand{\lim}{\varprojlim}

\newcommand{\elim}{\textsc{elim}}

\newcommand{\Id}{\mathrm{Id}}

\newcommand{\types}{\vdash}

\usepackage{amssymb,stmaryrd,amsthm,amsmath,bussproofs}
\usepackage[latin1]{inputenc}
\usepackage{graphics}
\usepackage{float}

\newcommand{\myemph}[1]{\textbf{#1}}    

\newcommand{\type}{\texttt{type}}       

\newcommand{\app}{\texttt{app}}

\newcommand{\id}[1]{\texttt{Id}_{#1}} 
\newcommand{\idn}[2]{\texttt{Id}_{#1}^{#2}}
\newcommand{\judge}[3][]{#2\;\vdash_{#1}\;#3}

\newtheorem{theorem}{Theorem}[section]

\theoremstyle{definition}

\theoremstyle{remark}

\floatstyle{ruled}
\newfloat{mytable}{htbp}{lot}[section]
\floatname{mytable}{Figure}

\begin{document}

\title{Type Theory and Homotopy}
\author{Steve Awodey}
\date{} 
\thanks{
Thanks to Pieter Hofstra, Peter Lumsdaine, and Michael Warren for their contributions to this article, and to Per Martin-L\"of and Erik Palmgren for supporting this work over many years.}



\maketitle

\section{Introduction}

\noindent The purpose of this informal survey article is to introduce the reader to a new and surprising connection between Geometry, Algebra, and Logic, which has recently come to light in the form of an interpretation of the constructive type theory of Per Martin-L\"of into homotopy theory, resulting in new examples of certain algebraic structures which are important in topology.  This connection was discovered quite recently, and various aspects of it are now under active investigation by several researchers.
 (See \cite{Awodey:HMIT,AHW:MLC,Warren:PhD, Berg:TWOG,Gambino:ITWFS,Garner:TDMTT,benno-richard,Lumsdaine:WOITT,vdBG:TSMIT,VVnote}.)%

\subsection{Type theory}

Martin-L\"of type theory is a formal system originally intended to provide a rigorous framework for constructive mathematics \cite{MartinLof:ITT,MartinLof:ITT72,MartinLof:ITT84}.  It is an extension of the typed $\lambda$-calculus admitting dependent types and terms. Under the Curry-Howard correspondence \cite{Howard:FTNC}, one identifies types with propositions, and terms with proofs;
viewed thus, the system is at least as strong as second-order logic, and it is known to interpret constructive set theory \cite{Aczel:stengthML}.
  Indeed, Martin-L\"of type theory has been used successfully to formalize large parts of constructive mathematics, such as the theory of generalized recursive definitions \cite{Nordstrom:PMLTT,MartinLof:CMCP}.
 Moreover, it is also employed extensively as a framework for the development of high-level programming languages,  
in virtue of its combination of expressive strength and desirable proof-theoretic properties \cite{Nordstrom:PMLTT,Streicher:STT}.

In addition to simple types $A, B, \dots$ and their terms $ x:A \vdash b(x) : B$, the theory also has dependent types $x:A\vdash B(x)$, which are regarded as indexed families of types.  There are simple type forming operations $A\times B$ and $A\rightarrow B$, as well as operations on dependent types, including in particular the sum $\sum_{x:A}B(x)$ and product $\prod_{x:A}B(x)$ types (see the appendix for details).  The Curry-Howard interpretation of the operations $A\times B$ and $A\rightarrow B$ is as propositional conjunction and implication, of course; the dependent types  $x:A\vdash B(x)$ are predicates, or more generally, relations, $$x_1:A_1, \dots, x_n:A_n \vdash R(x_1,\dots,x_n),$$ and the sum $\sum$ and product $\prod$ operations are the existential $\exists$ and universal $\forall$ quantifiers, respectively. 

It is now natural to further extend the type theory with a primitive equality relation, corresponding to the equality formulas of first-order logic.  Specifically, given two terms $a,b$\/ of the same type $A$, one can form a new \myemph{identity type} $\id{A}(a,b)$, representing the proposition that $a$\/ and $b$\/ are equal; a term of this type thus represents a proof of the proposition that $a$ equals $b$.  In the intensional version of the theory, with which we are concerned here, one thus has two different notions of equality: \myemph{propositional equality} is the notion represented by the identity types, in that two terms are propositionally equal just if their identity type $\id{A}(a,b)$ is inhabited by a term.   By contrast, \myemph{definitional equality} is a primitive relation on terms and is not represented by a type; it behaves much like equality between terms in the simply-typed lambda-calculus, or any conventional equational theory. 

If the terms $a$\/ and $b$\/ are definitially equal, then (since they can be freely substituted for each other) they are also propositionally equal; but the converse is generally not true in the intensional version of the theory (the rules for identity types are given in the appendix).  In the extensional theory, by contrast, the two notions of equality are forced by an additional rule to coincide. As a consequence, the extensional version of the theory is essentially a dependent type theory with a standard, extensional equality relation. As is well-known, however, the price one pays for this simplification is a loss of desirable proof-theoretic properties, such as strong normalization and decidable type checking and equality of terms \cite{Streicher:IIITT,Streicher:STT,Hofmann:ECITT}.%

In the intensional theory, each type $A$ is thus endowed by the identity types $\id{A}(a,b)$ with a non-trivial structure.  Indeed, this structure was observed by Hofmann and Streicher in \cite{Hofmann:GITT} to satisfy conditions analogous to the familiar laws for groupoids.\footnote{%
A \myemph{groupoid} is like a group, but with a partially-defined composition operation.  Precisely, a groupoid can be defined as a category in which every arrow has an inverse.  A group is thus a groupoid with only one object.  Groupoids arise in topology as generalized fundamental groups, not tied to a choice of basepoint (see below).}
  Specifically, the posited reflexivity of propositional equality produces identity proofs $\mathtt{r}(a):\id{A}(a,a)$\/ for any term $a:A$, playing the role of a unit arrow $1_a$ for $a$; and when $f:\id{A}(a,b)$\/ is an identity proof, then (corresponding to the symmetry of identity) there also exists a proof $f^{-1}:\id{A}(b,a)$, to be thought of as the inverse of $f$; finally, when $f:\id{A}(a,b)$\/ and $g:\id{A}(b,c)$\/ are identity proofs, then (corresponding to transitivity) there is a new proof $g \circ f:\id{A}(a,c)$, thought of as the composite of $f$\/ and $g$. Moreover, this structure on each type $A$ can be shown to satisfy the usual groupoid laws, but significantly,  only \myemph{up to propositional equality}.  We shall return to this point below.

The constructive character, computational tractability, and proof-theoretic clarity of the type theory are owed in part to this rather subtle treatment of equality between terms, which itself is expressible within the theory using the identity types $\id{A}(a,b)$.  Unlike extensional equality, which is computationally intractable, the expressibility of intensional equality within the theory leads to a system that is both powerful and expressive while retaining its important computational character.  The cost of intensionality, however, has long been the resulting difficulty of finding a natural, conventional semantic interpretation.  (See \cite{Hofmann:SSDT,cartmell:generalised-algebraic-theories,Hofmann:OITTLCCC,Dybjer:ITT} for previous semantics).

The new approach presented here constructs a bridge from constructive type theory to algebraic topology, exploiting both the axiomatic approach to homotopy of Quillen model categories, as well as the related algebraic methods involving (weak) higher-dimensional groupoids.  This at once provides \myemph{two} new domains of interpretation for type theory.  In doing so, it also permits logical methods to be combined with the traditional algebraic and topological approaches to homotopy theory, opening up a range of possible new applications of type theory in homotopy and higher-dimensional algebra.  It also allows the importation into homotopy theory of computational tools based on the type theory, such as the computer proof assistants Coq and Agda (cf.~\cite{Thery:C}).

\subsection{Homotopy theory}

In homotopy theory one is concerned with spaces and continuous mappings up to homotopy; a \myemph{homotopy} between continuous maps $f,g:X \to Y$\/ is a continuous map $\vartheta:X \times [0,1] \to Y$\/ satisfying $\vartheta(x,0)=f(x)$\/ and $\vartheta(x,1)=g(x)$. Such a homotopy $\vartheta$ can be thought of as a ``continuous deformation" of $f$ into $g$. Two spaces are said to be homotopy-equivalent if there are continuous maps going back and forth, the composites of which are homotopical to the respective identity mappings.  Such spaces may be thought of as differing only by a continuous deformation.  Algebraic invariants, such as homology or the fundamental group, are homotopy-invariant, in that any spaces that are homotopy-equivalent must have the same invariants. 

It is natural to also consider homotopies between homotopies, referred to as \myemph{higher homotopies}.  When we consider a space $X$, a distinguished point $p\in X$, and the paths in $X$ beginning and ending at $p$, and identify such paths up to homotopy, the result is the \myemph{fundamental group} $\pi(X,p)$ of the space at the point.  Pursuing an idea of Grothendieck's \cite{Grothendieck:PS}, modern homotopy theory generalizes this classical construction in several directions: first, we remove the dependence on the base-point $p$ by considering the \myemph{fundamental groupoid} $\pi(X)$, consisting of all points and all paths up to homotopy.  Next, rather than identifying homotopic paths, we can consider the homotopies between paths as distinct, new objects of a higher dimension (just as the paths themselves are homotopies between points).  Continuing in this way, we obtain a structure consisting of the points of $X$, the paths in $X$, the homotopies between paths, the higher homotopies between homotopies, and so on for even higher homotopies.  The resulting structure $\pi_\infty(X)$ is called the \myemph{fundamental weak $\infty$-groupoid of $X$}.  Such higher-dimensional algebraic structures now play a central role in homotopy theory (see e.g.~\cite{Kapranov:OGHT}); they capture much more of the homotopical information of a space than does the fundamental group $\pi(X,p)$, or the groupoid $\pi(X)=\pi_1(X)$, which is a quotient of $\pi_\infty(X)$ by collapsing the higher homotopies.  As discussed in subsection \ref{subsec:groupoid} below, it has recently been shown that such higher-dimensional groupoids also arise naturally in intensional type theory.

Another central concept in modern homotopy theory is that of a \myemph{Quillen model structure}, which captures axiomatically some of the essential features of homotopy of topological spaces, enabling one to ``do homotopy" in different mathematical settings, and to express the fact that two settings carry the same homotopical information.  
Quillen \cite{Quillen:HA} introduced model categories as an abstract framework for homotopy theory 
which would apply to a wide range of mathematical settings. 
Such a structure consists of the specification of three classes of maps (the fibrations, weak equivalences, and cofibrations) satisfying certain conditions typical of the leading topological examples.    
The resulting framework of axiomatic homotopy theory allows the development of the main lines of classical homotopy theory (fundamental groups, homotopies of maps, strong and weak equivalence, homotopy limits, etc.) independently of any one specific setting.  Thus, for instance, it is also applicable not only in spaces and simplicial sets, but also in new settings, as in the work of Voevodsky on the homotopy theory of 
schemes \cite{Morel:A1HTS}, or that of Joyal \cite{Joyal:QCKC,Joyal:QC} and Lurie \cite{Lurie:HTT} on quasicategories.  In the work under consideration here (subsection \ref{subsec:qmint}), it is shown that Martin-L\"of type theory can be interpreted in any model category. This allows the use of type theory to reason formally and systematically about homotopy theory.

%
\section{The homotopy interpretation}


\subsection{Background}

Among the most thorough, recent treatments of the \myemph{extensional} type theory are the two papers \cite{MP:WFTC,MP:TTCST} %
by Moerdijk and Palmgren from 2000 and 2002.  The authors also announced a projected third paper devoted to the intensional theory, which never appeared.  Their intention was presumably to make use of higher categories and, perhaps, Quillen model structures.  No preliminary results were stated, but see \cite{Palmgren:GLCC}.

In 2006, Vladimir Voevodsky gave a series of lectures at Stanford University entitled ``Homotopy lambda-calculus", in which an interpretation of intensional type theory into simplicial sets was proposed (see \cite{VVnote}).  At the same time, and independently, the author and his doctoral student Michael Warren established the interpretation of intensional type theory in Quillen model structures, following a suggestion of Moerdijk. 

All of these approaches derive from the pioneering work of Hoffmann and Streicher~\cite{Hofmann:GITT}, which we now summarize. 
\subsection{Groupoid semantics}

A model of type theory is \emph{extensional} if the following reflection rule is 
satisfied: 

\begin{prooftree}\label{rule:ref}
  \AxiomC{${ p : \id{A}(a, b)}$}
  \RightLabel{Id-reflection}
  \UnaryInfC{${ a = b : A}$}
\end{prooftree}
I.e., the identity type $\id{A}(a, b)$ in extensional models captures no more information than 
whether or not the terms $a$ and $b$ are definitionally equal. Although type checking is decidable in the intensional theory, it fails to be so in the extensional theory obtained by adding Id-reflection as a rule governing identity types. This fact is the principal motivation for studying intensional rather than extensional type theories (cf.~\cite{Streicher:STT} for a discussion of the difference between the intensional and extensional forms of the theory). 
A good notion of a model for the extensional theory is due to Seely~\cite{Seely:LCCCTT}, who showed that one can interpret type dependency in locally cartesian closed categories in a very natural way. (There are certain coherence issues, prompting a later refinement by Hofmann \cite{Hofmann:SSDT}, but this need not concern us here.)  Of course, intensional type theory can also be interpreted this way, but then the interpretation of the identity types necessarily becomes trivial in the above sense

The first natural, non-trivial semantics for intensional type theory were developed by Hoffmann and Streicher~\cite{Hofmann:GITT} using \myemph{groupoids}, which are categories in which every arrow is an iso.   The category of groupoids is not locally cartesian closed \cite{Palmgren:GLCC}, and the model employs certain fibrations (equivalently, groupoid-valued functors) to model type dependency.  Intuitively, the identity type over a groupoid $G$ is interpreted as the groupoid $G^{\rightarrow}$ of arrows in $G$, so that an identity proof $f:\id{A}(a,b)$\/ becomes an arrow $f:a\rightarrow b$ in $G$.  The interpretation no longer validates extensionality, since there can be different elements $a, b$ related by non-identity arrows $f:a\rightarrow b$.  Indeed, there may be many different such arrows $f,g: a\rightrightarrows b$\,; however---unlike in the type theory---these cannot in turn be further related by identity terms of higher type $\vartheta:\id{\id{A}}(f,g)$, since a (conventional) groupoid generally has no such higher-dimensional structure.  Thus the groupoid semantics validates a certain truncation principle, stating that all higher identity types are trivial---a form of extensionality one dimension up. In particular, the groupoid laws for the identity types are strictly satisfied in these models, rather than holding only up to propositional equality.

This situation suggests the use of the higher-dimensional analogues of groupoids, as formulated in homotopy theory, in order to provide models admitting non-trivial higher identity types.   Such higher groupoids occur naturally as the (higher) fundamental groupoids of spaces (as discussed above).   A step in this direction was made by Garner~\cite{Garner:TDMTT}, who uses a 2-dimensional notion of fibration to model intensional type theory in 2-groupoids, and shows that when various truncation axioms are added, the resulting theory is sound and complete with respect to this semantics.  In his dissertation \cite{Warren:PhD}, Warren showed that infinite-dimensional groupoids also give rise to models, which validate no such additional truncation axioms (see also \cite{Warren:OG}).  Such models do, however, satisfy type-theoretically unprovable strictness conditions such as the associativity of composition.  It seems clear that one will ultimately need to use \emph{weak} infinite dimensional groupoids in order to faithfully model the full intensional type theory (see subsection \ref{subsec:groupoid} below).

\subsection{Homotopical models of type theory}\label{subsec:qmint}

Groupoids and their homomorphisms arise in homotopy theory as a ``model" (i.e.\ a representation) of topological spaces with homotopy classes of continuous maps.  There are other models as well, such as simplicial sets.  The idea of a Quillen model structure (cf. \cite{Quillen:HA,Bousfield:CFSC}) is to axiomatize the common features of these different models of homotopy, allowing one to develop the theory in an abstract general setting, and to compare different particular settings.  

This axiomatic framework also provides a convenient way of specifying a general semantics for intensional type theory, not tied to a particular choice of groupoids, 2-groupoids, $\infty$-groupoids, simplicial sets, etc., or even spaces themselves. The basic result in this connection states that it is possible to model the intensional type theory in any Quillen model category \cite{Awodey:HMIT} (see also \cite{Warren:PhD}).  The idea is that a type is interpreted as an abstract ``space" $X$ and a term $x:X \vdash a(x):A$ as a continuous function $a : X\to A$.  Thus e.g.\ a closed term $a:A$ is a point $a$ of $A$, an identity term $p : \id{A}(a,b)$ is then a path $p:a \leadsto b$ in $A$ (a homotopy between points!). A ``higher" identity term  $\vartheta:\id{\id{A}(a,b)}(p,q)$ is a homotopy between the paths $p$ and $q$, and so on for even higher identity terms and higher homotopies.  In this interpretation, one uses abstract ``fibrations" to interpret dependent types, and abstract ``path spaces'' to model identity types, recovering the groupoid model and its relatives as special cases. 

In \cite{Gambino:ITWFS} 
it was then shown that the type theory itself carries a natural homotopy structure (i.e.\ a weak factorization system), so that the theory is not only sound, but also logically complete with respect to such abstract homotopical semantics.  While some ``coherence" issues regarding the strictness of the interpretation remain to be worked out (again, see \cite{Warren:PhD}, as well as \cite{vdBG:TSMIT}),  together these results clearly establish not only the viability of the homotopical interpretation as a semantics for type theory, but also the possibility of using type theory to reason in Quillen model structures.  That is to say, they suggest that intensional type theory can be seen as a ``logic of homotopy theory".

In order to describe the interpretation in somewhat more detail, we first recall a few standard definitions.
In any category $\C$, given maps $f : A\to B$ and $g : C\to D$, we write $f\pitchfork g$ 
to indicate that $f$ has the \emph{left-lifting property} (LLP) with respect to $g$: for any commutative square 
\[  
\xymatrix{
   A   
    \ar[r]^{h}
    \ar[d]_{f}
    &
   C
    \ar[d]^{g}
    \\
    B 
    \ar@{.>}[ru]|{j}
    \ar[r]_{i}
    &
    D
  }
\]
there exists a diagonal map  $j: B\to C$ such that $j\circ f = h$ and $g\circ j = i$. If $\M$ is any 
collection of maps in $\C$, we denote by $\leftexp{\pitchfork}{\M}$ the collection of maps in 
$\C$ having the LLP with respect to all maps in $\M$. The collection of maps $\M^{\pitchfork}$ is defined similarly. 
A \emph{weak factorization system} $(\L, \R)$ in a category $\C$ consists of two collections $\L$ (the ``left-class") and $\R$ (the ``right-class") of maps in $\C$ such that:\
\begin{enumerate} 
\item Every map $f : A\to B$ has a factorization as $f=p\circ i$, where $i \in \L$ and $p \in\R$. 
\[
\xymatrix{
A 
\ar[r]^i
\ar[rd]_f
&
C 
\ar[d]^p
\\
& B ,
}
\] 

\item $\L= \leftexp{\pitchfork}{\R}$ and $\L^\pitchfork =\R$. 
\end{enumerate}

A \emph{(closed) model category} \cite{Quillen:HA} is a bicomplete category $\C$ equipped with subcategories 
$\F$ (fibrations), $\sfC$ (cofibrations) and $\W$ (weak equivalences), satisfying the following two 
conditions: (1) Given any maps $g\circ f = h$, if any two of $f, g, h$ are weak equivalences, then so is the third; (2)
both $(\sfC, \F \cap \W)$ and $(\sfC \cap \W, \F)$ are weak factorization systems. 
A map $f$ in a model category is a \emph{trivial cofibration} if it is both a cofibration and a weak 
equivalence. Dually, a \emph{trivial fibration} is a map which is both a fibration and a weak equivalence. An object $A$ is said to be \emph{fibrant} if the canonical map $A\to1$ is a fibration. Dually, $A$ is \emph{cofibrant} if $0\to A$ is a cofibration. 

Examples of model categories include the following:
\begin{enumerate} 
\item The category $\mathsf{Top}$ of topological spaces, with fibrations the Serre fibrations, weak 
equivalences the weak homotopy equivalences, and cofibrations those maps which 
have the LLP with respect to trivial fibrations. The cofibrant objects in this model 
structure are retracts of spaces constructed, like CW-complexes, by attaching cells. 

\item The category $\mathsf{SSet}$ of simplicial sets, with cofibrations the monomorphisms, fibrations the Kan fibrations, and weak equivalences the weak homotopy equivalences. The fibrant objects for this model structure are the Kan complexes. 

\item The category $\mathsf{Gpd}$ of (small) groupoids, with cofibrations the homomorphisms that are
injective on objects, fibrations the Grothendieck fibrations, and weak equivalences the categorical equivalences. 
Here all objects are both fibrant and cofibrant. 
\end{enumerate}
See e.g.\ \cite{Dwyer:HTMC,Hovey:MC} for further examples and details.

Finally, recall that in any model category $\C$, a \emph{(very good) path object}  $A^I$ for an object  $A$ consists of a factorization 
\begin{equation}\label{diag:factor_Delta}
\xymatrix{
A 
\ar[r]^r
\ar[rd]_\Delta
&
A^I 
\ar[d]^p
\\
& A\times A ,
}
\end{equation}
of the diagonal map $\Delta : A \to A \times A$ as a trivial cofibration $r$ followed by a fibration 
$p$ (see \cite{Hovey:MC}).  Paradigm examples of path objects are given by exponentiation by a suitable ``unit interval"
$I$ in either $\mathsf{Gpd}$ or, when the object $A$ is a Kan complex, in $\mathsf{SSet}$.  
In e.g.\ the former case, $G^I$ is just the ``arrow groupoid'' $G^{\rightarrow}$, consisting of all arrows in the groupoid $G$.
Path objects always exist, but are not uniquely determined.  In many examples, however, they can be chosen functorially. 

We can now describe the homotopy interpretation of type theory more precisely. Whereas the idea of the Curry-Howard correspondence is often summarized by the slogan ``Propositions as Types", the idea underlying the homotopy interpretation is instead
``Fibrations as Types".  In classical topology, and in most model categories, a fibration $p : E\to X$ can be thought of as a family of objects $E_x$ varying continuously in a parameter $x\in X$. (The path-lifting property of a topological fibration describes how to get from one fiber $E_x = p^{-1}(x)$ to another $E_y$ along a path $f:x\leadsto y$).  This notion gives the interpretation of type dependency. 
Specifically, assume that  $\C$ is a finitely complete category with (at least) a weak factorization system $(\L, \R)$. Because most interesting examples arise from model categories, we refer to maps in $\L$ as trivial cofibrations and those in $\R$ as fibrations. A judgement $\vdash A : \type$ is then interpreted as a fibrant object $A$ of $\C$. Similarly, a dependent type $x : A\vdash B(x) : \type$ is interpreted as a fibration $p: B\to A$. Terms $ x: A\vdash b(x): B(x)$ in context are interpreted as sections $b:A\to B$ of $p: B\to A$, i.e. $p\circ b = 1_A$.  Thinking of fibrant objects as types and fibrations as dependent types, the natural interpretation of the identity type $\id{A}(a, b)$ should then be as the  \emph{fibration of paths} in $A$ from $a$ to $b$, so that the type $x,y:A\vdash \id{A}(x,y)$ should be the ``fibration of all paths in $A$". That is, it should be a path object for $A$. 

\begin{theorem}[\cite{Awodey:HMIT}]
Let $\C$ be a finitely complete category with a weak factorization system and a functorial choice of \emph{stable} path objects $A^I$:
i.e., given any fibration $A\to X$ and any map $f : Y\to X$, the evident comparison map is an isomorphism,
\[
f^*(A^I)\cong f^*(A)^I.
\]
Then $\C$ is a model of Martin-L\"of type theory with identity types. 
\end{theorem}

The proof exhibits the close connection between type theory and axiomatic reasoning in this setting: We verify the rules for the identity types (see the Appendix). Given a fibrant object $A$, the judgement $x, y : A\vdash \id{A}(x, y)$ is interpreted as the path object fibration $p : A^I \to A\times A$, see \eqref{diag:factor_Delta}. Because $p$ is then a fibration, the formation rule
\[
x, y : A\vdash \id{A}(x, y):\type
\]
is satisfied. Similarly, the introduction rule
\[
x:A \vdash \mathtt{r}(x):\id{A}(x,x)
\]
is valid because the interpretation $r : A\to A^I$ is a section of $p$ over $\Delta:A\to A\times A$. For the elimination and conversion rules, assume that the following premisses are given 
\begin{align*}
x : A, y : A, z : \id{A}(x, y) &\vdash D(x, y, z ) : \type , \\
x : A &\vdash d(x) : D(x, x, \mathtt{r}(x)) .
\end{align*}
We have, therefore, a fibration $q : D\to A^I$ together with a map $d : A\to D$ such that 
$q\circ d = r$. This data yields the following (outer) commutative square: 
\[  
\xymatrix{
   A   
    \ar[r]^{d}
    \ar[d]_{r}
    &
   D
    \ar[d]^{q}
    \\
    A^I 
    \ar@{.>}[ru]|{j}
    \ar[r]_{1}
    &
    A^I
  }
\]
 Because $q$ is a fibration and  $r$ is, by definition, a trivial cofibration, there exists a diagonal filler $j$, which we choose
as the interpretation of the term: 
\[
x, y : A, z : \id{A}(x, y) \vdash \mathtt{J}(d, x, y, z ) : D(x, y, z ). 
\]
Commutativity of the bottom triangle is precisely this conclusion of the elimination rule, and commutativity of the top triangle is the required conversion rule:
\[
 x : A \vdash \mathtt{J}(d, x, x, \mathtt{r}(x)) = d(x) : D(x, x, \mathtt{r}(x)).
\] 

Examples of categories satisfying the hypotheses of this theorem include groupoids, simplicial sets, and 
many simplicial model categories \cite{Quillen:HA} (including, e.g., simplicial sheaves and presheaves).  There is a question of selecting the diagonal fillers $j$ as interpretations of the $\mathtt{J}$-terms in a ``coherent way", i.e.\ respecting substitutions of terms for variables.  Some solutions to this problem are discussed in \cite{Awodey:HMIT,Warren:PhD,Garner:CGNWFS}.  One neat solution is implicit in the recent work of Riehl \cite{RiehlAMS} on ``algebraic" Quillen model structures.  A systematic investigation of the issue of coherence, along with several examples of coherent models derived from homotopy theory, can be found in the recent work \cite{vdBG:TSMIT} of van den Berg and Garner.

\subsection{Higher algebraic structures}\label{subsec:groupoid}

Given the essential soundness and completeness of type theory with respect to the homotopical interpretation, we may further ask how \emph{expressive} the logical system is, as a language for homotopy theory?  From this point of view, we think of the types in the intensional theory as spaces,  the terms of the type $A$\/ are the points of the ``space" $A$, the identity type $\id{A}(a,b)$\/ represents the collection of paths from $a$\/ to $b$, and the higher identities are homotopies between paths,  homotopies between homotopies of paths, etc., and we ask what homotopically relevant facts, properties, and structures are logically expressible. The topological fact that paths and homotopies do not form a groupoid, but only a groupoid up to homotopy, is of course reminiscent of the logical fact that the identity types only satisfy the groupoid laws up to propositional equality. This apparent \emph{analogy} between homotopy theory and type theory can now be made precise, and indeed can be recognized as one and the same fact, resting entirely on the homotopical interpretation of the logic.  The fundamental weak $\omega$-groupoid of a space is namely a construction entirely within the logical system --- it belongs, as it were, to the logic of homotopy theory, as we now proceed to explain.

\subsubsection{Weak $\omega$-groupoids}

It has recently been shown by Peter Lumsdaine~\cite{Lumsdaine:WOITT} and, independently, Benno van den Berg and Richard Garner~\cite{Berg:TWOG,benno:talk}, that the tower of identity types over any fixed base type $A$\/ in the type theory bears an infinite dimensional algebraic structure of exactly the kind arising in homotopy theory, called a weak $\omega$-groupoid (\cite{Kapranov:OGHT,Leinster:survey,cheng:duals-give-inverses,Brown:FGTG}).

In somewhat more detail,  in the globular approach to higher groupoids \cite{Leinster:HOHC,Batanin:MGCNETWNC}, a weak $\omega$-groupoid has objects (``0-cells''), arrows (``1-cells'') between objects, 2-cells between 1-cells, and so on, with various composition operations and laws depending on the kind of groupoid in question (strict or weak, $n$- or $\omega$-, etc.).  
We first require the notion of a globular set, which may be thought of as an ``infinite-dimensional" graph.  Specifically, a \emph{globular set} (\cite{Batanin:MGCNETWNC,Street:PTGS}) is a presheaf on the category $\G$ generated by arrows
$$ 0 \two^{s_0}_{t_0} 1 \two^{s_1}_{t_1} 2 \two \ldots $$
subject to the equations $ss = ts$, $st = tt$.  
More concretely, a globular set $\A$ has a set $A_n$ of ``$n$-cells'' for each $n \in \N$, and each $(n+1)$-cell $x$ has parallel source and target $n$-cells $s(x)$, $t(x)$.  (Cells $x,y$ of dimension $>0$ are \emph{parallel} if $s(x) = s(y)$ and $t(x) = t(y)$; all $0$-cells are considered parallel.)  

\begin{figure}[htp]
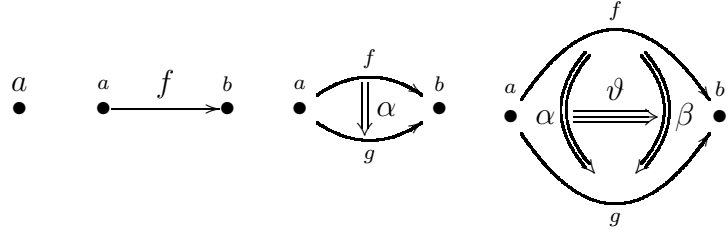

$$
\begin{array}{c}
\begin{array}{cccc}
\ \xy
(0,0)*{\bullet};
(0,80)*{a};
\endxy \quad
&
\ \xy
(0,0)*{\bullet}="a";
(0,80)*{\scriptstyle a};
(400,0)*{\bullet}="b";
(400,80)*{\scriptstyle b};
{\ar "a";"b"};
(200,80)*{f};
\endxy \ 
&
\ \xy
(0,0)*+{\bullet}="a";
(0,80)*{\scriptstyle a};
(450,0)*+{\bullet}="b";
(450,80)*{\scriptstyle b};
{\ar@/^1pc/^{f} "a";"b"};
{\ar@/_1pc/_{g} "a";"b"};
{\ar@{=>} (210,85)*{};(210,-85)*{}};
(280,0)*{\alpha};
\endxy \ 
&
\xy 0;/r.22pc/: 
(0,15)*{}; 
(0,-15)*{}; 
(0,8)*{}="A"; 
(0,-8)*{}="B"; 
{\ar@{=>}@/_.75pc/ "A"+(-4,1) ; "B"+(-3,0)}; 
(-10,0)*{\alpha};
{\ar@{=}@/_.75pc/ "A"+(-4,1) ; "B"+(-4,1)}; 
{\ar@{=>}@/^.75pc/ "A"+(4,1) ; "B"+(3,0)}; 
(10,0)*{\beta};
{\ar@{=}@/^.75pc/ "A"+(4,1) ; "B"+(4,1)}; 
{\ar@3{->} (-6,0)*{} ; (6,0)*+{}}; 
(0,4)*{\vartheta};
(-15,0)*+{\bullet}="1"; 
(-15,4)*{\scriptstyle a};
(15,0)*+{\bullet}="2"; 
(15,4)*{\scriptstyle b};
{\ar@/^2.75pc/^{f} "1";"2"}; 
{\ar@/_2.75pc/_{g} "1";"2"}; 
\endxy
\end{array} 
\end{array}
$$
\caption{Some cells in dimensions 0--3. \label{figure:assoc-laws}} 
\end{figure}

For example, given a type $A$ in a type theory $\T$, the terms of types $$A, \id{A}, \id{{\id{A}}}, \dots,$$ together with the evident indexing projections, e.g.\  $s(p) = a$ and $t(p) = b$ for $p:\id{A}(a,b)$, form a globular set $\hat{A}$.

A strict $\omega$-groupoid is an infine-dimensional groupoid satisfying, in all dimensions, associativity, unit, and inverse laws given by equations between certain cells.   Such a groupoid has an underlying globular set consisting of cells of each dimension, and any globular set $\A$ generates a free strict $\omega$-groupoid $F(\A)$---just as any set generates a free group, and any graph, a free groupoid.  The cells of $F(\A)$ are free (strictly associative) pastings-together of cells from $\A$ and their formal duals, including degenerate pastings from the identity cells of $F(\A)$.  In a \emph{strict} $\omega$-groupoid, cells can be composed along a common boundary in any lower dimension, and the composition satisfies various associativity, unit, and interchange laws, captured by the generalized associativity law: each labelled pasting diagram has a unique composite.

In a \emph{weak} $\omega$-groupoid, by contrast, we do not expect strict associativity, and so we may have multiple composition maps for each shape of pasting diagram; but we do demand that these composites agree \emph{up to cells of the next dimension}, and that these associativity cells satisfy coherence laws of their own, and so on.

Now, this is exactly the situation we find in intensional type theory.  For instance, even in constructing a term witnessing the transitivity of identity, one finds that there is no single canonical candidate.  Specifically, as a composition for the pasting diagram
$$ \xymatrix{ \cdot \ar[r] & \cdot \ar[r] & \cdot }$$
or more concretely, a term $c$ such that 
$$x,y,z:X, p:\Id(x,y), q:\Id(y,z) \types c(q,p): \Id(x,z),$$
there are the two equally natural terms $c_l$, $c_r$ obtained by applying ($\Id$-\elim) to $p$ and $q$ respectively.  These are not definitionally equal, but are propositionally equal, i.e.\ equal up to a 2-cell, for there is a term $e$ with
$$x,y,z:X, p:\Id(x,y), q:\Id(y,z) \types e(q,p): \Id(c_l(q,p),c_r(q,p)).$$
Indeed, we have the following:

\begin{theorem}[\cite{Lumsdaine:WOITT,Berg:TWOG}]  Let $A$ be any type in a system  $\T$ of intensional Martin-L\"of type theory.  Then the globular set $\hat{A}$ of terms of type $$A, \id{A}, \id{\id{A}}, \ldots$$ carries a natural weak $\omega$-groupoid structure.
\end{theorem}

It is now quite natural to ask what special properties this particular $\omega$-groupoid has in virtue of its type-theoretic construction.  In light of related syntactic constructions of other types of free algebras, a reasonable conjecture is that it is the \myemph{free weak $\omega$-groupoid}, up to a suitable notion of equivalence.  We return to this question below.

\subsubsection{Weak $n$-groupoids}

A further step in exploring the connection between type theory and homotopy is to investigate the relationship between type theoretic ``truncation" (i.e.\ higher-dimensional extentionality principles) and topological ``truncation" of the higher fundamental groups. Spaces for which the homotopy type is already completely determined by the fundamental groupoid are called \myemph{homotopy 1-types}, or simply 1-types \cite{Baues:HT}. More generally, one has n-types, which are thought of as spaces which have no homotopical information above dimension n. One of the goals of homotopy theory is to obtain good models of homotopy n-types.  For example, the category of groupoids is Quillen equivalent to the category of 1-types; in this precise sense, groupoids are said to model homotopy 1-types.  A famous conjecture of Grothendieck's is that (arbitrary) homotopy types are modeled by weak $\infty$-groupoids (see e.g.~\cite{Batanin:MGCNETWNC} for a precise statement).

Recent work \cite{AHW:MLC} by the author, Pieter Hofstra, and Michael Warren has shown that the 1-truncation of the intensional theory, arrived at by adding the analogue of the Id-reflection rule for all terms of identity type, generates a Quillen model structure on a category of structured graphs that is Quillen equivalent to that of groupoids.  In a precise sense, the truncated system of 1-dimensional type theory thus models the homotopy 1-types.


In a bit more detail, for every globular set $\A$ one can define a system of type theory  $\T(\A)$, the basic terms of which are the elements of the various $A_n$, typed as terms of the corresponding identity type determined by the globular structure: e.g.\ $a\in A_n$ is a basic term of type $\id{A}(s(a), t(a))$, where $s,t :A_n \rightrightarrows A_{n-1}$  are the source and target maps, at dimension n, of $\A$.  Since we know from the result of Lumsdaine et al. \cite{Lumsdaine:WOITT,Berg:TWOG},
just reviewed, that for any type $X$, the underlying globular set of terms of the various identity types $X, \id{X}, \id{\id{X}}, \dots$ gives rise to a weak $\omega$-groupoid, we can infer that in particular the globular set of terms over the ground type $A_0$ in the theory $\T(\A)$ form such a groupoid, \myemph{generated type-theoretically} from the arbitrary globular set $\A$.  Let us call this weak $\omega$-groupoid $G_{\omega}(\A)$, the \myemph{type-theoretically free} weak $\omega$-groupoid generated by $\A$.  This construction is investigated in depth in \cite{AHW:MLC},  where certain groupoids of this kind are termed \myemph{Martin-L\"of complexes} (technically, these are the algebras for the globular monad just described).  

It is clearly of interest to investigate the relationship between this type-theoretic construction of higher groupoids and both the algebraically free higher groupoids, on the one hand, and the higher group\-oids arising from spaces as fundamental groupoids, on the other.  As a first step, one can consider the 1-dimensional truncation of the above construction, and the resulting (1-) groupoid $G_{1}(\A)$.  For that case, the following result relating $G_{1}(\A)$ to the usual, algebraically free groupoid is established in the work cited:

\begin{theorem}[\cite{AHW:MLC}] The type-theoretically free groupoid is equivalent to the algebraically free groupoid.
\end{theorem}

Furthermore, it is shown that the 1-truncated Martin-L\"of complexes admit a Quillen model structure equivalent to that of (1-) groupoids.  The following then results from known facts from homotopy theory:

\begin{theorem}[\cite{AHW:MLC}] The 1-truncated Martin-L\"of complexes classify homotopy 1-types.
\end{theorem}

Obviously, one should now proceed to higher groupoids and the corresponding type theories truncated at higher dimensions.

\section{Conclusion: The logic of homotopy}

The application of logic in geometry and topology via categorical algebra has a precedent in the development of topos theory.  Invented by Grothendieck as an abstract framework for sheaf cohomology, the notion of a topos was soon discovered to have a logical interpretation, admitting the use of logical methods into topology (see e.g. \cite{JoyalTierney:EGTG} for just one of many examples).  Equally important was the resulting flow of geometric and topological ideas and methods into logic, e.g.\ sheaf-theoretic independence proofs, topological semantics for many non-classical systems, and an abstract treatment of realizability (see the encyclopedic work \cite{Johnstone:E1}).

An important and lively research program in current homotopy theory is the pursuit (again following Grothendieck \cite{Grothendieck:PS}) of a general concept of ``stack," subsuming sheaves of homotopy types, higher groupoids, quasi-categories, and the like.  Two important works in this area have just appeared (Lurie, \emph{Higher Topos Theory} \cite{Lurie:HTT}; Joyal, \emph{Theory of Quasi-Categories} \cite{Joyal:QC}).  It may be said, somewhat roughly, that the notion of a ``higher-dimensional topos" is to homotopy what that of a topos is to topology (as in \cite{Joyal:SSCS}).  This concept also has a clear categorical-algebraic component via Grothendieck's ``homotopy hypothesis", which states that $n$-groupoids are combinatorial models for homotopy $n$-types, and $\infty$-groupoids are models for arbitrary homotopy types of spaces.  Still missing from the recent development of higher-dimensional toposes, however, is a logical aspect analogous to that of (1-dimensional) topos theory.   The research surveyed here suggests that such a logic  is already available in intensional type theory.  The homotopy interpretation of Martin-L\"of type theory into Quillen model categories, and the related results on type-theoretic constructions of higher groupoids, are analogous to the basic results interpreting \emph{extensional} type theory and higher-order logic in (1-) toposes.  They clearly indicate that the logic of higher toposes---i.e., the logic of homotopy---is, rather remarkably, a form of intensional type theory.

\appendix
\section{Rules of type theory}

\noindent This appendix recalls (some of) the rules of intensional Martin-L\"of type theory.
See \cite{MartinLof:ITT84,Nordstrom:PMLTT,Jacobs:CLTT} for detailed presentations.


\subsection*{Judgement forms}

There are four basic  forms of judgement:
\begin{align*}
A:\type &\qquad a:A\\
a=b:A &\qquad A = B:\type
\end{align*}
Each form can occur also with free variables: e.g.\ if $A$ is a type, then
\[
x:A \vdash B(x):\type
\]
is called a \emph{dependent type}, regarded as an $A$-indexed family of types. 
The part $x:A$ to the left of the turnstile $\vdash$ is called the \emph{context} of the judgement.   
More generally, a list of variable declarations $x_1:A_1, x_2:A_2,\ldots, x_n:A_n$ is a context whenever the judgements $A_1:\type$ and $$x_1 : A_1, . . . , x_m : A_m \vdash A_{m+1}:\type$$ are derivable, for $1 \leq m < n$. 
Given such a context $\Gamma$, the judgement $\Gamma\vdash A:\type$ means that $A$ is a type (in context $\Gamma$), while $\Gamma\vdash a:A$ indicates that $a$ is a term of type $A$ (in context $\Gamma$); the equality judgements have their usual meaning.





\subsection*{Formation rules}
Given an $A$-indexed family of types $B(x)$, the dependent sum $\sum_{x:A}.B(x)$ and product $\prod_{x:A}.B(x)$ can be formed. The identity type introduces a new dependent type $\id{A}$ for any type $A$.

\begin{prooftree}
  \AxiomC{$\judge{x:A}{B(x):\type}$}
  \RightLabel{$\prod$ formation}
  \UnaryInfC{$\prod_{x:A}B(x):\type$}
\end{prooftree}
\smallskip
\begin{prooftree}
  \AxiomC{$\judge{x:A}{B(x):\type}$}
  \RightLabel{$\sum$ formation}
  \UnaryInfC{$\sum_{x:A}B(x):\type$}
\end{prooftree}
\smallskip
\begin{prooftree}
  \AxiomC{$A:\type$}
  \RightLabel{$\id{}$ formation}
  \UnaryInfC{$\judge{x:A, y:A}{\idn{A}{}(x,y):\type}$}
\end{prooftree}
\medskip

\noindent Under the Curry-Howard correspondence, sums correspond to existential quantifiers, products to universal quantifiers, and identity types to equations. The behavior of each of these types is specified by introduction, elimination and conversion rules.

\subsection*{Rules for dependent products} 

\begin{prooftree}
  \AxiomC{$\judge{x:A}{f(x):B(x)}$}
  \RightLabel{$\prod$ introduction}
  \UnaryInfC{$\lambda x.f(x):\prod_{x:A}B(x)$}
\end{prooftree}
\smallskip
\begin{prooftree}
  \AxiomC{$a:A$}
  \AxiomC{$f:\prod_{x:A}B(x)$}
  \RightLabel{$\prod$ elimination}
  \BinaryInfC{$\app(f,a):B(a).$}
\end{prooftree}
\smallskip
\begin{prooftree}
  \AxiomC{$a:A$}
  \AxiomC{$\judge{x:A}{f(x):B(x)}$}
  \RightLabel{$\prod$ conversion}
  \BinaryInfC{$\app\bigl(\lambda x.f(x),a\bigr) \;=\; f(a):B(a)$}
\end{prooftree}
\medskip

\noindent The introduction rule states that for every family of terms $f(x) : B(x)$ there is a term $\lambda x.f(x)$ of type $\prod_{x:A}B(x)$. The elimination rule corresponds to the application of a term $f$ of the indexed product to $a : A$. 
Finally, the conversion rule for states that the application term $\app(-, a)$ behaves correctly when applied to a term of the form $\lambda x.f(x)$.

\subsection*{Rules for dependent sums} 

\begin{prooftree}
  \AxiomC{$a:A$}
  \AxiomC{$b:B(a)$}
  \RightLabel{$\sum$ introduction}
  \BinaryInfC{$\langle a,b\rangle:\sum_{x:A}B(x)$}
\end{prooftree}
\smallskip
\begin{prooftree}
  \AxiomC{${p:\sum_{x:A}B(x)}$}
  \AxiomC{$\judge{x:A, y:B(x)}{c(x,y):C(\langle x,y\rangle)}$}
  \RightLabel{$\sum$ elimination}
  \BinaryInfC{$\sigma(c,p):C(p)$}
\end{prooftree}
\smallskip
\begin{prooftree}
  \AxiomC{$a:A$}
  \AxiomC{$b:B(a)$}
  \AxiomC{$\judge{x:A,y:B(x)}{c(x,y):C(\langle x,y\rangle)}$}
  \RightLabel{$\sum$ conversion}
  \TrinaryInfC{$\sigma(c,\langle a,b\rangle) \;=\; c(a,b):C(\langle a,b\rangle)$}
\end{prooftree}
\medskip

\noindent The variables $x:A, y:B(a)$ are bound in the the notation $\sigma(c,p)$.  

Note that when $A$ and $B$ are types in the same context, the usual product $A \times B$ and function $A\rightarrow B$ types from the simply typed $\lambda$-calculus are recovered as $\sum_{x:A}B$ and $\prod_{x:A}B$, respectively. 

\subsection*{Rules for identity types}
%
\begin{prooftree}
  \AxiomC{$a:A$}
  \RightLabel{$\id{}$ introduction}
  \UnaryInfC{$\mathtt{r}(a):\idn{A}{}(a,a)$}
\end{prooftree}
\medskip
\begin{prooftree}
 \AxiomC{$c:\idn{A}{}(a,b)$}
  \AxiomC{$\judge{x:A,y:A,z:\idn{A}{}(x,y)}{B(x,y,z):\type}$}
  \noLine
  \UnaryInfC{$\judge{x:A}{d(x):B\bigl(x,x,\mathtt{r}(x)\bigr)}$}
  \RightLabel{$\id{}$ elimination}
  \BinaryInfC{$\mathtt{J}(d,a,b,c):B(a,b,c)$}
\end{prooftree}
\smallskip
\begin{prooftree}
  \AxiomC{$a:A$}
  \RightLabel{$\id{}$ conversion}
  \UnaryInfC{$\mathtt{J}\bigl(d,a,a,\mathtt{r}(a)\bigr)\;=\;d(a):B\bigl(a,a,\mathtt{r}(a)\bigr)$}
\end{prooftree}
\medskip

\noindent The introduction rule provides a witness $\mathtt{r}(a)$ that $a$ is identical to itself, called the \emph{reflexivity term}. 
The distinctive elimination rule can be recognized as a form of Leibniz's law.
The variable $x:A$ is bound in the the notation $\mathtt{J}(d,a,b,c)$.
\bibliographystyle{alpha}
\bibliography{htt}

\end{document}